\documentclass[10pt]{amsart}
\usepackage{subfigure}
\usepackage{graphicx}
\usepackage{enumerate}
\usepackage{bbm}
\usepackage[
colorlinks=true, urlcolor=blue]{hyperref}

\newtheorem{thm}{Theorem}
\newtheorem{thm*}{Theorem}
\newtheorem{lem}[thm]{Lemma}
\newtheorem{cor}[thm]{Corollary}
\newtheorem{conj}[thm]{Conjecture}
\newtheorem{prop}[thm]{Proposition}

\theoremstyle{remark}

\theoremstyle{definition}
\newtheorem{deef}[thm]{Definition}

\newcommand{\R}{\mathbbm{R}}
\newcommand{\rd}{\mathrm{d}}

\usepackage{color}

\title{The Willmore energy and the magnitude of Euclidean domains}
\author{Heiko Gimperlein, Magnus Goffeng}
\address{Heiko Gimperlein\newline
\indent Leopold-Franzens-Universit\"{a}t Innsbruck\newline
\indent Technikerstraße 13 \newline
\indent 6020 Innsbruck\newline
\indent Austria \newline\newline
\indent Magnus Goffeng,\newline
\indent Centre for Mathematical Sciences\newline 
\indent University of Lund\newline 
\indent Box 118, SE-221 00 Lund\newline 
\indent Sweden\newline}
\email{heiko.gimperlein@uibk.ac.at, magnus.goffeng@math.lth.se}
\begin{document}
\maketitle

\begin{abstract}
We study the geometric significance of Leinster's notion of magnitude for a compact metric space. For a smooth, compact domain $X$ in an odd-dimensional Euclidean space, we show that the asymptotic expansion of the function $\mathcal{M}_X(R) = \mathrm{Mag}(R\cdot X)$ at $R = \infty$ determines the Willmore energy of the boundary $\partial X$. This disproves the Leinster-Willerton conjecture for a compact convex body in odd dimensions.
\end{abstract}

\section*{Introduction}

The notion of magnitude was introduced by Leinster \cite{leinsterold,leinster} as an extension of the Euler characteristic to (finite) enriched categories. Magnitude has been shown to unify notions of ``size'' like the cardinality of a set, the length of an interval or the Euler characteristic of a triangulated manifold, and it even relates to measures of the diversity of a biological system. See \cite{leinmeck} for an overview.

Viewing a metric space as a category enriched over $[0,\infty)$, Leinster and Willerton proposed and studied the magnitude of metric spaces \cite{leinster, leinwill}: If $(X,d)$ is a finite metric space, a weight function is a function $w : X \to \R$ which satisfies $\sum_{y \in X} \mathrm{e}^{-\rd(x,y)}w(y) =1$ for all $x \in X$. Given a weight function $w$, we define the magnitude of $X$ as $\mathrm{Mag}(X) := \sum_{x \in X} w(x)$; this definition is independent of the choice of weight function. Beyond finite metric spaces, the magnitude of a compact, positive definite metric space $(X,\rd)$ was made rigorous by Meckes \cite{meckes}: 
$$\mathrm{Mag}(X) := \sup\{\mathrm{Mag}(\Xi) : \Xi \subset X\ \ \text{finite}\}\ .$$
Instead of the magnitude of an individual space $(X,d)$, it proves fruitful to study the magnitude function $\mathcal{M}_X(R) := \mathrm{Mag}(X, R\cdot d)$ for $R>0$. 

Compact convex subsets $X \subset \mathbb{R}^n$ provide a key example, surveyed in \cite{leinmeck}. Motivated by properties of the Euler characteristic and computer calculations, Leinster and Willerton \cite{leinwill} conjectured a surprising relation to the intrinsic volumes $V_i(X)$, which would shed light on the geometric content of the magnitude function:
\begin{equation}
\label{leinwillor}
\mathcal{M}_X(R)=\sum_{k=0}^n\frac{1}{k!\omega_k}V_k(X)R^k+o(1),\quad\mbox{as $R\to \infty$}.
\end{equation}
Here, $\omega_k$ is the volume of the $k$-dimensional unit ball. This asymptotic expansion resembles the well-known expansion of the heat trace, with leading terms $V_n(X)=\textnormal{vol}_n(X)$, $V_{n-1}(X)=\textnormal{vol}_{n-1}(\partial X)$ \cite{gilkey}. The expansion coefficients for the heat trace, however, are not proportional to $V_k(X)$ for $k\leq n-2$.

The conjectured behavior \eqref{leinwillor} was disproved by Barcel\'o and Carbery \cite{barcarbs} for the unit ball $B_5 \subset \mathbb{R}^5$. They explicitly computed the rational function $\mathcal{M}_{B_5}$ and observed numerical disagreement of the coefficients of $R^k$. Their results were extended to balls in odd dimensions in \cite{will}. 

In spite of this negative result, the authors were able to prove a variant of \eqref{leinwillor}, with modified prefactors, which confirmed the close relation between magnitude and intrinsic volumes \cite{gimpgoff}: When $n=2m-1$ is odd and $X\subseteq \mathbb{R}^n$ is a compact domain with smooth boundary, there are coefficients $(c_j(X))_{j\in \mathbb{N}}$ such that 
$$\mathcal{M}_X(R)=\sum_{j=0}^\infty\frac{c_j(X)}{n!\omega_n}R^{n-j}+O(R^{-\infty}),\quad\mbox{as $R\to \infty$},$$
where $$ c_0(X)=\textnormal{vol}_n(X),\  c_1(X)=m\textnormal{vol}_{n-1}(\partial X),\  c_2(X)=\frac{m^2}{2}\ (n-1)\int_{\partial X} H\, \mathrm{d} S\ .$$
Here, $H$ denotes the mean curvature of $\partial X$. Each coefficient $c_{j}$ is an integral over $\partial X$ computable from the second fundamental form of $\partial X$ and its covariant derivatives. For $j=0,1,2$ and $X$ convex, the coefficient $c_j$ is proportional to the intrinsic volume $V_{n-j}(X)$, for $j=0,1,2$. This proves that the Leinster-Willerton conjecture holds for modified universal coefficients up to $O(R^{n-3})$. 

The following variant of the Leinster-Willerton conjecture therefore  remained plausible. It would confirm the relation between magnitude and intrinsic volumes and, in particular, show that $c_n$ is proportional to the Euler characteristic $V_0$: 
\begin{conj}
\label{wekawleie}
For $n>0$, there are universal constants $\gamma_{0,n},\gamma_{1,n},\ldots, \gamma_{n,n}$ such that for any compact convex subset $X\subseteq \mathbb{R}^n$, $\mathcal{M}_X(R)=\sum_{k=0}^n\gamma_{k,n}V_k(X)R^k+o(1)$, as $R\to \infty$.
\end{conj}
In this paper we prove that Conjecture \ref{wekawleie} fails in all odd dimensions $n\geq 3$ and find unexpected geometric content in $c_3$. While the conjecture holds true for the terms of order $R^n$, $R^{n-1}$ and $R^{n-2}$, the $R^{n-3}$-term is not proportional to an intrinsic volume: 
\begin{thm}
\label{computc3}
Assume that $n\geq 3$ is odd and that $X\subseteq \mathbb{R}^n$ is a compact domain with smooth boundary. Then there is a dimensional constant $\lambda_n\neq 0$ such that 
$$c_3(X)=\lambda_n\mathcal{W}(\partial X),$$
where 
$\mathcal{W}(\partial X):=\int_{\partial X} H^2\mathrm{d}S$
is the Willmore energy of the boundary of the hypersurface $\partial X$.
\end{thm}
Building on \cite{gimpgoff}, the proof reformulates the magnitude function in terms of an elliptic boundary value problem of order $n+1$ in $\R^n\setminus X$, which is then studied using methods from semiclassical analysis. See Proposition \ref{magformula} and Equation \eqref{ckdirneu} below.\\

To see that Theorem \ref{computc3} disproves Conjecture \ref{wekawleie} in the fourth term, we observe that the Willmore energy is not an intrinsic volume: The only intrinsic volume with the same scaling property as the Willmore energy is $V_{n-3}$. For instance, if $n=3$ then $V_{n-3}$ is the Euler characteristic while $\int_{\partial X}H^2\mathrm{d}S$ can be non-zero even when $\partial X$ has vanishing Euler characteristic (e.g. for a torus). In general dimension, for $a>0$ the solid ellipsoid 
$$X_a:=\left\{(x',x_n)\in \mathbb{R}^{n-1}\times \mathbb{R}: |x'|^2+\frac{|x_n|^2}{a^2}\leq 1\right\},$$
satisfies that $\mathcal{W}(\partial X_a)\to \infty$ as $a\to 0$. On the other hand, Hausdorff continuity of intrinsic volumes shows that $V_{n-3}(X_a)$ converges to a finite number, namely the $n-3$:rd intrinsic volume of the $n-1$-dimensional unit ball. Therefore Theorem \ref{computc3} implies the following.

\begin{cor}
Assume that $n \geq 3$ is odd and that $X\subseteq \mathbb{R}^n$ is a compact convex domain with smooth boundary. There are universal constants $\gamma_{n-2,n},\gamma_{n-1,n}, \gamma_{n,n}$ such that
$$\mathcal{M}_X(R)=\sum_{k=n-2}^n\gamma_{k,n}V_k(X)R^k+O(R^{n-3}),\quad\mbox{as $R\to \infty$}.$$
However, there is no constant $\gamma_{n-3,n}$ such that 
$\mathcal{M}_X(R)=\sum_{k=n-3}^n\gamma_{k,n}V_k(X)R^k+O(R^{n-4})$ as $R\to \infty$.
In particular, the Leinster-Willerton conjecture fails even with modified universal coefficients.
\end{cor}

\subsection*{Acknowledgements} {The authors thank the anonymous referee for their feedback, which helped improve the paper.} MG was supported by the Swedish Research Council Grant VR 2018-0350. We thank Tom Leinster for comments on an earlier draft.

\section*{Background and notation}

We assume that $X\subseteq \R^n$ is a compact domain with $C^\infty$-boundary, where $n=2m-1$ odd. Denote by $\Omega:=\R^n\setminus X$ \emph{the exterior domain}. We use the Sobolev spaces $H^s(\R^n):=(1-\Delta)^{-s/2}L^2(\R^n)$ of exponent $s\geq 0$. Here, the Laplacian $\Delta$ is given by $\Delta=\sum_{j=1}^n \frac{\partial^2}{\partial x_j^2}$. The spaces $H^s(X)$ and $H^s(\Omega)$ are defined using restrictions. The Sobolev spaces $H^s(\partial X)$ can be defined using local charts or as $(1-\Delta_{\partial X})^{-s/2}L^2(\partial X)$.

We use $\partial_\nu$ to denote the Neumann trace of a function $u$ in $\Omega$. The operator $\partial_\nu$ extends to a continuous operator $H^s(\Omega)\to H^{s-3/2}(\partial X)$ for $s>3/2$. Similarly, $\gamma_0:H^s(\Omega)\to H^{s-1/2}(\partial X)$ denotes the trace operator defined for $s>1/2$. 

For $R>0$ we shall need the operators 
$$\mathcal{D}^j_R:=\begin{cases}
\partial_\nu \circ(R^2-\Delta)^{(j-1)/2}, \; &\mbox{when $j$ is odd},\\
\gamma_0\circ(R^2-\Delta)^{j/2}, \; &\mbox{when $j$ is even}.
\end{cases}$$
By the trace theorem, $\mathcal{D}^j_R$ is continuous as an operator $\mathcal{D}^j_R:H^s(\Omega)\to H^{s-j-1/2}(\partial X)$ for $s>j+1/2$. 

We recall a key observation from \cite{barcarbs}, in the reformulation presented in \cite{gimpgoff}:

\begin{prop}\cite[Proposition 9]{gimpgoff}
\label{magformula}
Suppose that $h_R\in H^{2m}(\Omega)$ is the unique weak solution to the boundary value problem
\begin{align*}
\begin{cases}
(R^2-\Delta)^mh_R&=0\quad\mbox{in $\Omega$}\\
\vspace{-3mm}\\
\mathcal{D}^j_R h_R&=
\begin{cases} 
R^j, \;& j\mbox{   even}\\
0, \;& j\mbox{   odd}.
\end{cases}
\;\;, \; j=0,...,m-1.
\end{cases}
\end{align*}
Then the following identity holds 
$$\mathcal{M}_X(R)=\frac{\textnormal{vol}_n(X)}{n!\omega_n}R^n-\frac{1}{n!\omega_n}\sum_{\frac{m}{2}<j\leq m} R^{n-2j}\int_{\partial X}\mathcal{D}^{2j-1}_R h_R\,\rd S.$$
\end{prop}

The operators $\mathcal{D}^j_R$ define a matrix-valued Dirichlet-Neumann operator $\Lambda(R):\mathcal{H}_+\to \mathcal{H}_-$ in the Hilbert space
$$\mathcal{H}:=\underbrace{\bigoplus_{j=0}^{m-1}H^{2m-j-1/2}(\partial X)}_{\mathcal{H}_+}\oplus \underbrace{\bigoplus_{j=m}^{n}H^{2m-j-1/2}(\partial X)}_{\mathcal{H}_-}$$
as follows:  $\Lambda(R)(u_j)_{j=0}^{m-1}:=(\mathcal{D}^j_Ru)_{j=m}^n,$ where $u\in H^{2m}(\Omega)$ is the unique weak solution to 
\begin{equation}
\label{thebvpgener}
\begin{cases}
(R^2-\Delta)^mu&=0\quad\mbox{in $\Omega$}\\
\mathcal{D}^j_R u&=u_j
\;\;, \; j=0,...,m-1.
\end{cases}
\end{equation}
The operator $\Lambda(R)$ is a parameter-dependent pseudodifferential operator on $\partial X$. {  The parameter $R$ enters like an additional co-variable, which allows us to compute the asymptotics of $\mathcal{M}_X$ from Proposition \ref{magformula}. For the convenience of the reader we recall the salient features of the parameter-dependent pseudodifferential calculus, see for instance \cite{grubb77,GrubbGreenBook,shubinbook} for further details. We restrict to parameters $R\in \R_+=(0,\infty)$.

\begin{deef}
\label{defineidnomom}
A parameter-dependent pseudodifferential operator $A$ of order $s$ on $\R^n$ is an operator on the Schwartz space of the form 
\begin{equation}
\label{knokjnad}
Af(x):=(2\pi)^{-n}\int_{\R^n}\int_{\R^n} a(x,\xi,R)\mathrm{e}^{i(x-y)\xi}f(y)\mathrm{d}y\mathrm{d}\xi, \quad f\in \mathcal{S}(\R^n),
\end{equation}
where the full symbol $a$ admits a polyhomogeneous expansion of order $s$ in $(\xi,R)$. That is, for $k\in \mathbb{N}$ there are functions $a_{s-k}\in C^\infty(\R^n\times \R^n\times \R_+)$ with 
$$a_{s-k}(x,t\xi,tR)=t^{s-k}a_{s-k}(x,\xi,R), \quad\mbox{for $t\geq 1$, $\|(\xi,R)\|\geq 1$},$$ 
and $a$ can be written as an asymptotic sum
$$a\sim \sum_{k=0}^\infty a_{s-k}.$$
We call $a_s$ the principal symbol of $A$. If $a_s(x,\xi,R)$ is invertible for every $(x,\xi,R)\in \R^n\times\R^n\times \R_+$, we say that $A$ is elliptic with parameter.
\end{deef}

Definition \ref{defineidnomom} on $\R^n$ extends by standard techniques, using coordinate charts, to define a pseudodifferential operator and its full symbol on a compact manifold, see for instance \cite{gimpgoff,grubb77,GrubbGreenBook,shubinbook}. The use of the parameter-dependent calculus is crucial to the work \cite{gimpgoff} and the computations in this paper, including formulas for the symbol of a product of two pseudodifferential operators and the parametrix construction. In particular, if $A$ is elliptic with parameter of order $s$ on a compact manifold, it has a parametrix with parameter $B$ of order $-s$. The full symbol expansion $b\sim \sum_{j=0}^\infty b_{-s-j}$ can be explicitly computed: The principal symbol is given by $b_{-s}(x,\xi,R)=a_{s}(x,\xi,R)^{-1}$, and for $j>0$ the following inductive formula holds,
\begin{equation}
\label{parametricpknoknad}
b_{-s-j}(x,\xi,R)=-a_{s}(x,\xi,R)^{-1}\!\!\!\sum_{\substack{k+l+|\alpha|=j,\\ l<j}}\frac{i^{|\alpha|}}{\alpha!} \partial_\xi^\alpha a_{s-k}(x,\xi,R)\partial_x^\alpha b_{-s-l}(x,\xi,R).
\end{equation}
The computation proving Equation \eqref{parametricpknoknad} follows from \cite[Section 5.5]{shubinbook}.

For $R \to \infty$ the parameter-dependent calculus further allows to compute expectation values of the form $\int_M A(1)\ \mathrm{d}x$ in terms of the symbol:

\begin{lem}
\label{asymptoomadoma}
Suppose that $A:C^\infty(M)\to C^\infty(M)$ is a parameter-dependent pseudodifferential operator of order $s$ acting on a compact manifold $M$ equipped with a volume density. Then there is an asymptotic expansion 
$$\int_M A(1)\ \mathrm{d}x=\sum_{k=0}^\infty a_k R^{s-k} +O(R^{-\infty}),$$
where the coefficients $a_k$ are computed as follows: Expand the full symbol of $A$ into terms homogeneous in $(\xi,R)$ as $\sigma_A(x,\xi,R)\sim \sum_{k=0}^\infty \sigma_{s-k}(A)(x,\xi,R)$ and set 
$$a_k:=\int_M \sigma_{s-k}(A)(x,0,1)\ \mathrm{d}x.$$
\end{lem}

For the proof of Lemma \ref{asymptoomadoma} we refer the reader to \cite[Lemma 20]{gimpgoff} or \cite[Lemma 2.24]{gimpgofflouc}, but let us outline the main idea. The claimed asymptotics of Lemma \ref{asymptoomadoma} is coordinate invariant because $\int_M A(1)\ \mathrm{d}x$ is coordinate invariant. It therefore suffices to compute the asymptotics for an operator $A$ on $\R^n$ as in Equation \eqref{knokjnad}, assuming $a$ is compactly supported in the $x$-variable. In this case, $A(1)=a(x,0,R)$, so that for $R\geq 1$
\begin{align*}
\int_{\R^n} A(1)\ \mathrm{d}x=&\int_{\R^n}a(x,0,R)\mathrm{d}x=\sum_{k=0}^\infty \int_M \sigma_{s-k}(A)(x,0,R)\ \mathrm{d}x+O(R^{-\infty})=\\
=&\sum_{k=0}^\infty \int_M \sigma_{s-k}(A)(x,0,1)\ \mathrm{d}xR^{s-k}+O(R^{-\infty})=\sum_{k=0}^\infty a_k R^{s-k} +O(R^{-\infty}).
\end{align*}
The reader should note that the integrands $a_{s-k}(x,0,R)=a_{s-k}(x,0,1)R^{m-k}$ are well defined because each $a_{s-k}$ is homogeneous in $(\xi,R)$, and not only in $\xi$.

From Proposition \ref{magformula} and Lemma \ref{asymptoomadoma} }we deduce a formula for the expansion coefficients $c_k$:
\begin{equation}
\label{ckdirneu}
c_k(X):=-\sum_{\frac{m}{2}<j\leq m}\sum_{0\leq l<m/2}\int_{\partial X}\sigma_{2j-2l-k}(\Lambda_{2j-1,2l})(x,0,1)\,\mathrm{d} S,
\end{equation}
for $k>0$ where $\Lambda=(\Lambda_{j+m,l})_{j,l=0}^{m-1}$ and $\sigma_{2j-2l-k}(\Lambda_{2j-1,2l})$ the homogeneous part of order $2j-2l-k$ in its symbol (with parameter). See \cite[Proposition 20]{gimpgoff}.

{ The full symbol of the parameter-dependent operator $\Lambda$ can be computed by adapting} standard techniques in semiclassical analysis \cite{GrubbGreenBook}. The operator $\Lambda$ is first computed using boundary layer potentials. To define these, we consider the function 
$$K(R;z):=\frac{\kappa_n}{R} \mathrm{e}^{-R|z|}, \quad z\in \R^n.$$
The constant $\kappa_n>0$ is chosen such that
$$(R^2-\Delta)^mK=\delta_0$$  in the sense of distributions on $\R^n$. For $l=0,\ldots,n$, we define the functions
$$K_l(R;x,y):=(-1)^l\mathcal{D}^{n-l}_{R,y}K(R;x-y), \quad x\in \R^n, \; y\in \partial X.$$
Here $\mathcal{D}^l_{R,y}$ denotes $\mathcal{D}^l_R$ acting in the $y$-variable. We also consider the distributions
$$K_{j,k}(R;x,y):= \mathcal{D}^j_{R,x}K_{k}(R;x,y), \quad x\in \partial X.$$
Each $K_{j,k}$ defines a { parameter-dependent pseudodifferential operator} $A_{j,k}(R):C^\infty(\partial X)\to C^\infty(\partial X)$, 
\begin{align*}
A_{j,k}(R)f(x)&:=\int_{\partial X} K_{j,k}(R;x,y)f(y)\mathrm{d}S(y), \quad x\in \partial X.
\end{align*}
The integral defining $A_{j,k}(R)$ is understood in the sense of an exterior limit.
These operators combine into a $2m\times 2m$-matrix of operators $\mathbb{A}:=(A_{j,l})_{j,l=0}^n:\mathcal{H}\to \mathcal{H}$. It decomposes into matrix blocks
$$\mathbb{A}=\begin{pmatrix} \mathbb{A}_{++}& \mathbb{A}_{+-}\\ \mathbb{A}_{-+}& \mathbb{A}_{--}\end{pmatrix}:\begin{matrix} \mathcal{H}_+\\\oplus\\\mathcal{H}_-\end{matrix}\longrightarrow \begin{matrix} \mathcal{H}_+\\\oplus\\\mathcal{H}_-\end{matrix},$$
with $\mathbb{A}_{pq}:\mathcal{H}_q\to \mathcal{H}_p$ for $p,q\in \{+,-\}$. { By integrating by parts as in \cite[Proposition 12]{gimpgoff}, one can show that if $u$ solves Equation \eqref{thebvpgener} then 
$$u_+=\mathbb{A}_{++}u_++\mathbb{A}_{+-}u_-,$$
where $u_+:=(u_j)_{j=0}^{m-1}$ and $u_-:=(u_{m+j})_{j=0}^{m-1}$. Therefore, $(1-\mathbb{A}_{++})u_+=\mathbb{A}_{+-}u_-$ and} we can express the Dirichlet-Neumann operator $\Lambda$ in terms of layer potentials as 
\begin{equation}
\label{DNlayer}
\Lambda=\mathbb{B}(1-\mathbb{A}_{++}).
\end{equation}
Here $\mathbb{B}=(B_{j+m,l})_{j,l=0}^{m-1}$ denotes a parametrix (with parameter) of $\mathbb{A}_{+-}=(A_{j,l+m})_{j,l=0}^{m-1}$. See more in the proof of \cite[Theorem 18]{gimpgoff}. 

The proof of Theorem \ref{computc3} uses Equation \eqref{DNlayer} to compute components of the symbol of the Dirichlet-Neumann operator $\Lambda$. The formula for $c_3$ then follows from \eqref{ckdirneu}.

\section*{Proof of Theorem \ref{computc3}}

{ To prove Theorem \ref{computc3} we note that we by Equation \eqref{ckdirneu} only need to compute the third term $\sigma_{2j-2l-3}(\Lambda_{2j-1,2l})$ in the polyhomogeneous expansion
$$\sigma(\Lambda_{2j-1,2l})(x,\xi,R)\sim \sum_{k=0}^\infty \sigma_{2j-2l-1-k}(\Lambda_{2j-1,2l})(x,\xi,R),$$
in the range $\frac{m}{2}<j\leq m$, $0\leq l<m/2$. In fact, we only need to compute the evaluation $\sigma_{2j-2l-3}(\Lambda_{2j-1,2l})(x,0,1)$. Recall that we are using the parameter-dependent calculus, so that each $\sigma_{2j-2l-1-k}(\Lambda_{2j-1,2l})(x,\xi,R)$ is homogeneous of degree $-2j-2l-1-k$ in $(\xi,R)$. 

For the convenience of the reader, we change to the notation $(x',\xi',R)\in T^*\partial X\times \mathbb{R}_+$ for coordinates and cotangent variables on the boundary $\partial X$, as used in \cite{gimpgoff}. For an integer $k\in \mathbb{Z}$, we use the notation 
\begin{align*}
\sigma_{k}(\mathbb{A}_{++})&:=(\sigma_{j-l+k}(A_{j,l}))_{j,l=0}^{m-1},\\
\sigma_{k}(\mathbb{A}_{+-})&:=(\sigma_{j-l+k-m}(A_{j,l+m}))_{j,l=0}^{m-1}\quad\mbox{and}\\
 \sigma_{k}(\mathbb{B})&:=(\sigma_{j+m-l+k}(B_{j+m,l}))_{j,l=0}^{m-1}.
\end{align*}
Here we write $\sigma_{j-l+k}(A_{j,l})$ for the degree $j-l+k$ part of $a_{j,l}$ written as a symbol depending on the variable $(x',\xi',R)\in T^*\partial X\times \mathbb{R}_+$. The symbols $\sigma_{k}(\mathbb{A}_{++})$, $\sigma_{k}(\mathbb{A}_{+-})$ and $\sigma_{k}(\mathbb{B})$ relate to the (parameter-dependent) Douglis-Nirenberg calculus naturally appearing in the boundary reduction of boundary value problems \cite{gimpgoff, grubb77}.} The reader should note the difference with the expressions appearing just after \cite[Proposition 37]{gimpgoff} in that they are for symbols in the variables $(x',y',\xi',R)$. The process of going between these two symbol expressions is one of the difficulties in the computation ahead.

The reader can note that $\sigma_{0}(\mathbb{A}_{++})$, $\sigma_{0}(\mathbb{A}_{+-})$ and $\sigma_{0}(\mathbb{B})$ are the matrices of principal symbols of $\mathbb{A}_{++}$, $\mathbb{A}_{+-}$ and $\mathbb{B}$, respectively. In particular, 
$$\sigma_{0}(\mathbb{B})=\sigma_{0}(\mathbb{A}_{+-})^{-1}.$$
It follows from \cite[Theorem 12]{gimpgoff} that $\sigma_{0}(\mathbb{B})$ does not depend on $x'\in \partial X$. Define the symbol 
$$\mathbb{D}=(\delta_{j,k}(R^2+|\xi|^2)^{j/2})_{j,k=0}^n.$$
By { the computational result} \cite[Theorem 12]{gimpgoff}, there are constant $m\times m$-matrices $C_0$, $C_1$, $C_2$, $C_3$ such that 
\begin{align*}
\sigma_{0}(\mathbb{A}_{++})&=\mathbb{D}C_0\mathbb{D}^{-1},\quad
\sigma_{0}(\mathbb{A}_{+-})=\mathbb{D}C_1\mathbb{D}^{-1},\\
\sigma_{-1}(\mathbb{A}_{++})&=H\mathbb{D}C_2\mathbb{D}^{-1},\quad
\sigma_{-1}(\mathbb{A}_{+-})=H\mathbb{D}C_3\mathbb{D}^{-1},\quad\mbox{and}\\
 \sigma_{0}(\mathbb{B})&=\mathbb{D}C_1^{-1}\mathbb{D}^{-1},
\end{align*}
where $H$ denotes the mean curvature of $\partial X$ and we in each identity embed $m\times m$-matrices in a suitable fashion into $2m\times 2m$-matrices. 

From \cite[Lemma 22, part a]{gimpgoff} and the $x'$-independence of $\sigma_{0}(\mathbb{B})$ we can { from Equation \eqref{parametricpknoknad}} deduce that 
$$\sigma_{-1}(\mathbb{B})=-\sigma_{0}(\mathbb{B})\sigma_{-1}(\mathbb{A}_{+-})\sigma_{0}(\mathbb{B})=H\mathbb{D}C_1^{-1}C_3C_1^{-1}\mathbb{D}^{-1},$$
as well as
\begin{align*}
\sigma_{-2}(\mathbb{B})=-\sigma_0(\mathbb{B})\bigg(\sigma_{-2}(\mathbb{A}_{+-})+\sum_{j=1}^{n-1}\partial_{\xi_j}\sigma_{0}(\mathbb{A}_{+-})\sigma_0(\mathbb{B})\partial_{x_j}&\sigma_{-1}(\mathbb{A}_{+-})-\\
-\sigma_{-1}(\mathbb{A}_{+-})\sigma_0&(\mathbb{B})\sigma_{-1}(\mathbb{A}_{+-})\bigg)\sigma_0(\mathbb{B}).
\end{align*}

Using \cite[Lemma 22, part b]{gimpgoff}, we write 
\begin{align*}
\sigma_{-2}(\Lambda)=&\sigma_{-2}(\mathbb{B})(1-\sigma_0(\mathbb{A}_{++}))-\sigma_{-1}(\mathbb{B})\sigma_{-1}(\mathbb{A}_{++})-\sigma_0(\mathbb{B})\sigma_{-2}(\mathbb{A}_{++})+\\
&+i\sum_{j=1}^{n-1}\partial_{\xi_j}\sigma_{-1}(\mathbb{B})\partial_{x_j}\sigma_{-1}(\mathbb{A}_{++})=\\
=&-\sigma_0(\mathbb{B})\bigg(\sigma_{-2}(\mathbb{A}_{+-})+\sum_{j=1}^{n-1}\partial_{\xi_j}\sigma_{0}(\mathbb{A}_{+-})\sigma_0(\mathbb{B})\partial_{x_j}\sigma_{-1}(\mathbb{A}_{+-})-\\
&\qquad\qquad-\sigma_{-1}(\mathbb{A}_{+-})\sigma_0(\mathbb{B})\sigma_{-1}(\mathbb{A}_{+-})\bigg)\sigma_0(\mathbb{B})(1-\sigma_0(\mathbb{A}_{++}))+\\
&+\sigma_{0}(\mathbb{B})\sigma_{-1}(\mathbb{A}_{+-})\sigma_{0}(\mathbb{B})\sigma_{-1}(\mathbb{A}_{++}))-\sigma_0(\mathbb{B})\sigma_{-2}(\mathbb{A}_{++})-\\
&-i\sum_{j=1}^{n-1}\partial_{\xi_j}\left(\sigma_{0}(\mathbb{B})\sigma_{-1}(\mathbb{A}_{+-})\sigma_{0}(\mathbb{B})\right)\partial_{x_j}\sigma_{-1}(\mathbb{A}_{++})
\end{align*}
Since all $\sigma_0$-occurences only depend on $R^2+|\xi|^2$, all its $\xi$-derivatives will vanish at $\xi=0$, and therefore, 
\begin{align*}
\sigma_{-2}(\Lambda)&(x',0,R)=\\
=&\bigg[-\sigma_0(\mathbb{B})\bigg(\sigma_{-2}(\mathbb{A}_{+-})-\sigma_{-1}(\mathbb{A}_{+-})\sigma_0(\mathbb{B})\sigma_{-1}(\mathbb{A}_{+-})\bigg)\sigma_0(\mathbb{B})(1-\sigma_0(\mathbb{A}_{++}))+\\
&+\sigma_{0}(\mathbb{B})\sigma_{-1}(\mathbb{A}_{+-})\sigma_{0}(\mathbb{B})\sigma_{-1}(\mathbb{A}_{++})-\sigma_0(\mathbb{B})\sigma_{-2}(\mathbb{A}_{++})\bigg]_{\xi'=0}=\\
=&\bigg[-\mathbb{D}C_1^{-1}\mathbb{D}^{-1}\left(\sigma_{-2}(\mathbb{A}_{+-})\mathbb{D}C_1^{-1}(1-C_0)\mathbb{D}^{-1}+\sigma_{-2}(\mathbb{A}_{++})\right)+\\
&+H^2\mathbb{D}C_1^{-1}C_3C_1^{-1}C_3C_1^{-1}(1-C_0)\mathbb{D}^{-1}+H^2\mathbb{D}C_1^{-1}C_3C_1^{-1}C_2\mathbb{D}^{-1}\bigg]_{\xi'=0}.
\end{align*}

Assume for now that $\sigma_{-2}(\mathbb{A}_{+-})(x',0,R)=\sigma_{-2}(\mathbb{A}_{++})(x',0,R)=0$. Then this computation shows that  indeed, there are universal constants $(d_{j+m,l})_{j,l=0}^{m-1}$ (independent of $X$) such that for $\frac{m}{2}<j\leq m$ and $0\leq l<m/2$, 
$$\sigma_{2j-2l-2}(\Lambda_{2j-1,2l})(x,0,1)=d_{2j-1,2l}H(x)^2.$$
In particular, we have shown that for a dimensional constant $\lambda_n$, we have that 
$c_3(X)=\lambda_n\int_{\partial X}H^2\mathrm{d}S$.
It follows from \cite{will} that $\lambda_n\neq 0$ for $n\geq 3$ odd.

It remains to show that $\sigma_{-2}(\mathbb{A}_{+-})(x',0,R)=\sigma_{-2}(\mathbb{A}_{++})(x',0,R)=0$. Note that we do not claim that $\sigma_{-2}(\mathbb{A}_{+-})=\sigma_{-2}(\mathbb{A}_{++})=0$ just that when restricting to $\xi'=0$ the symbols vanish. { This last step in the proof relies on the technically involved computations in \cite[Appendix A.2]{gimpgoff} and the process of going from ``two-variable symbols'' $\tilde{a}(x,y,\xi,R)$ to ``one-variable symbols'' $a(x,\xi,R)$, see \cite[Theorem 7.13]{GrubbDistOp}.} We pick local coordinates at a point on $\partial X$. We can assume that this point is $0\in \mathbb{R}^n$ and that the coordinates are of the form $(x',S(x'))$, where $x'$ belongs to some neighborhood of $0\in \mathbb{R}^{n-1}$ and $S$ is a scalar function with $S(0)=0$ and $\nabla S(0)=0$. We can express $a_{jk}$ as 
\begin{align*}
a_{jk}(x',y',\xi',R)=&b_{0,m-p-q}(R^2+|\xi'|^2,S(x')-S(y')), \\ 
&\quad\quad\mbox{when  $j=2p, k=n-2q$}\\
a_{jk}(x',y',\xi',R)=&b_{1,m-p-q}(R^2+|\xi'|^2,S(x')-S(y'))+\\
&(\xi'\cdot \nabla S(x'))b_{0,m-p-q}(R^2+|\xi'|^2,S(x')-S(y')), \\ 
&\quad\quad\mbox{when  $j=2p+1, k=n-2q$}
\end{align*}
\begin{align*}
a_{jk}(x',y',\xi',R)=&b_{1,m-p-q}(R^2+|\xi'|^2,S(x')-S(y'))+\\
&(\xi'\cdot \nabla S(y'))b_{0,m-p-q}(R^2+|\xi'|^2,S(x')-S(y')),\\ 
&\quad\quad\mbox{when  $j=2p, k=n-2q-1$}\\
a_{jk}(x',y',\xi',R)=&b_{2,m-p-q}(R^2+|\xi'|^2,S(x')-S(y'))+\\
&((\xi'\cdot \nabla S(y'))+(\xi'\cdot \nabla S(x')))b_{1,m-p-q}(R^2+|\xi'|^2,S(x')-S(y'))+\\
&(\xi'\cdot \nabla S(x'))(\xi'\cdot \nabla S(y'))b_{0,m-p-q}(R^2+|\xi'|^2,S(x')-S(y')), \\ 
&\quad\quad\mbox{when  $j=2p+1, k=n-2q-1$},
\end{align*} 
where 
\[
b_{r,N}(u,z)=
\begin{cases} 
(-i\partial_z)^r(u-\partial_z^2)^{-N}\delta_{z=0}, \;& N\leq 0,\\
{}\\
(-i\partial_z)^r\sum_{k=0}^{N-1}\tilde{c}_{k,r,N} \frac{|z|^k\mathrm{e}^{-|z|\sqrt{u}}}{u^{N-(k+1)/2}}, \;& N> 0,
\end{cases}
\]
for some coefficients $\tilde{c}_{k,r,N}$. 

We need to verify that $\sigma_{j-k-2}(A_{j,k})(x',0,R)=0$ for any $j$ and $k$. The symbol $\sigma_{j-k-2}(A_{j,k})$ in $x'=0$ is { by \cite[Theorem 7.13]{GrubbDistOp}} given by the terms of order $j-k-2$ in the expression 
$$a_{jk}(0,0,\xi',R)-i\sum_{l=1}^{n-1}\frac{\partial^2 a_{jk}}{\partial \xi_l\partial y_l}(0,0,\xi',R)-\frac{1}{2}\sum_{l,s=1}^{n-1}\frac{\partial^4 a_{jk}}{\partial \xi_l\partial \xi_s\partial y_l\partial y_s}(0,0,\xi',R)$$
Recall that $S(0)=0$ and $\nabla S(0)=0$ so there are several terms vanishing when setting $x'=0$. Indeed, no term of order $j-k-2$ in $a_{jk}(0,0,\xi',R)$ is non-zero. All non-zero terms of order $j-k-2$ in $\sum_{l=1}^{n-1}\frac{\partial^2 a_{jk}}{\partial \xi_l\partial y_l}(0,0,\xi',R)$ are odd functions under the reflection $\xi'\mapsto -\xi'$, so they vanish when restricting to $\xi'=0$. Similar computations show that terms of order $j-k-2$ in $\frac{1}{2}\sum_{l,s=1}^{n-1}\frac{\partial^4 a_{jk}}{\partial \xi_l\partial \xi_s\partial y_l\partial y_s}$ all contains a factor of $\xi_l$ or $\xi_l\xi_s$ so they vanish when restricting to $\xi'=0$.

\end{document}